\documentclass[12pt, reqno]{amsart}
\usepackage{amsmath, amsthm, amscd, amsfonts, amssymb, graphicx, color}
\usepackage[bookmarksnumbered, colorlinks, plainpages]{hyperref}

\hypersetup{colorlinks=true,linkcolor=red, anchorcolor=green, citecolor=cyan, urlcolor=red, filecolor=magenta, pdftoolbar=true}
\newtheorem{theorem}{Theorem}[section]

\theoremstyle{definition}
\newtheorem{definition}[theorem]{Definition}

\theoremstyle{remark}
\newtheorem{remark}[theorem]{Remark}
\numberwithin{equation}{section}
\input{tcilatex}

\begin{document}
\title[Existence of solutions for a class of quasilinear systems]{\textbf{%
Conditions for the existence of positive radial solutions for a class of
quasilinear systems}}
\author[ D.-P. Covei]{\textsf{Dragos-Patru Covei}$^{1}$\textsf{\ }}
\address{ $^{1}${\small \textit{Department of Applied Mathematics}}\\
{\small \ \textit{The Bucharest University of Economic Studies }}\\
{\small \textit{Piata Romana, 1st district, postal code: 010374, postal
office: 22, Romania}}}
\email{\textcolor[rgb]{0.00,0.00,0.84}{coveidragos@yahoo.com}}
\date{Received: xxxxxx; Revised: yyyyyy; Accepted: zzzzzz. \\
\indent $^{*}$Corresponding author}

\begin{abstract}
By using a monotone iterative scheme and Arzela-Ascoli theorem, we show the
existence of positive radial solutions to the quasilinear systems%
\begin{equation*}
\left\{ 
\begin{array}{c}
\Delta _{\phi _{1}}u:=a_{1}(|x|)f_{1}(v),\quad x\in \mathbb{R}^{N}\text{ (}%
N\geq 3\text{)}, \\ 
\Delta _{\phi _{2}}v:=a_{2}(|x|)f_{2}(u),\quad x\in \mathbb{R}^{N}\text{ (}%
N\geq 3\text{)},%
\end{array}%
\right.
\end{equation*}%
under appropriate conditions on the functions $\phi _{1}$, $\phi _{2}$, the
weights $a_{1}$, $a_{2}$ and to the nonlinearities $f_{1},$ $f_{2}$.
Moreover, we obtain a number of qualitative results concerning the behavior
of solutions. We also point that the functions $\phi _{1}$ and $\phi _{2}$
includes special cases appearing in mathematical models in nonlinear
elasticity, plasticity, generalized Newtonian fluids, and in quantum physics.
\end{abstract}

\maketitle

\setcounter{page}{1}


\let\thefootnote\relax\footnote{%
Copyright 2016 by the JOURNAL.}

\section{Introduction and Statement of the Main Results}

In this paper, we are concerned with the existence of nonnegative solutions
for a quasilinear system of the type 
\begin{equation}
\left\{ 
\begin{array}{l}
\Delta _{\phi _{1}}u:=a_{1}(\left\vert x\right\vert )f_{1}(v),\text{ }x\in 
\mathbb{R}^{N}\text{ (}N\geq 3\text{)}, \\ 
\Delta _{\phi _{2}}v:=a_{2}(\left\vert x\right\vert )f_{2}(u),\text{ }x\in 
\mathbb{R}^{N}\text{ (}N\geq 3\text{)},%
\end{array}%
\right.  \label{11}
\end{equation}%
where $\Delta _{\phi _{i}}u$ $(i=1,2)$ stands for the $\phi _{i}$-Laplacian
operator defined as $\Delta _{\phi _{i}}u:=\func{div}(\phi _{i}(|\nabla
u|)\nabla u)$ and the $C^{1}$-functions $\phi _{1}$ and $\phi _{2}$ satisfy
throughout the paper the following conditions:

(O1)\quad $\phi _{i}$ $\in C^{1}\left( \left( 0,\infty \right) ,\left(
0,\infty \right) \right) $ and $\underset{t\rightarrow 0}{\lim }t\phi
_{i}\left( t\right) =0$;

(O2)\quad $t\phi _{i}(t)>0$ is strictly increasing for $t>0$;

(O3)\quad there exist positive constants $\underline{k}_{i},\overline{k}_{i}$%
, the continuous and increasing functions $\underline{\theta }_{i}$, $%
\overline{\theta }_{i}:\left[ 0,\infty \right) \rightarrow \left[ 0,\infty
\right) $ and the continuous functions $\underline{\psi }_{i}$, $\overline{%
\psi }_{i}:\left[ 0,\infty \right) \rightarrow \left[ 0,\infty \right) $
such that 
\begin{equation}
\underline{k}_{i}\underline{\theta }_{i}(s_{1})\underline{\psi }%
_{i}(s_{2})\leq h_{i}^{-1}(s_{1}s_{2})\leq \overline{k}_{i}\overline{\theta }%
_{i}(s_{1})\overline{\psi }_{i}(s_{2})\text{ for all }s_{1},s_{2}>0,
\label{ineq}
\end{equation}%
where $h_{i}^{-1}$ is the inverse function of $h_{i}(t)=t\phi _{i}(t)$ for $%
t>0$.

The motivation for the present work stems from recent investigations of the
authors \cite{DF}, \cite{L1,L2}. We give a quick review here. Lair, \cite{L1}
has considered entire large radial solutions for the elliptic system%
\begin{equation}
\left\{ 
\begin{array}{l}
\Delta u=a_{1}\left( \left\vert x\right\vert \right) v^{\alpha }, \\ 
\Delta v=a_{2}\left( \left\vert x\right\vert \right) u^{\beta }\text{, }x\in 
\mathbb{R}^{N}\text{ (}N\geq 3\text{),}%
\end{array}%
\right.  \label{lair}
\end{equation}%
where $0<\alpha \leq 1$, $0<\beta \leq 1$, $a_{1}$ and $a_{2}$ are
nonnegative continuous functions on $\mathbb{R}^{N}$, and he proved that a
necessary and sufficient condition for this system to have a nonnegative
entire large radial solution (i.e., a nonnegative spherically symmetric
solution $\left( u,v\right) $ on $\mathbb{R}^{N}$ that satisfies $\underset{%
\left\vert x\right\vert \rightarrow \infty }{\lim }u\left( x\right) =%
\underset{\left\vert x\right\vert \rightarrow \infty }{\lim }v\left(
x\right) =\infty $), is 
\begin{eqnarray}
\int_{0}^{\infty }ta_{1}\left( t\right) \left(
t^{2-N}\int_{0}^{t}s^{N-3}Q\left( s\right) ds\right) ^{\alpha }dt &=&\infty ,
\label{c1l} \\
\int_{0}^{\infty }ta_{2}\left( t\right) \left(
t^{2-N}\int_{0}^{t}s^{N-3}P\left( s\right) ds\right) ^{\beta }dt &=&\infty ,
\label{c2l}
\end{eqnarray}%
where $P\left( r\right) =\int_{0}^{r}\tau a_{1}\left( \tau \right) d\tau $
and $Q\left( r\right) =\int_{0}^{r}\tau a_{2}\left( \tau \right) d\tau $.

It is well known, see Yang \cite{DF}, that if $a:\left[ 0,\infty \right)
\rightarrow \left[ 0,\infty \right) $ is a spherically symmetric continuous
function and the nonlinearity $f:[0,\infty )\rightarrow \lbrack 0,\infty )$
is a continuous, increasing function with $f\left( 0\right) \geq 0$ and $%
f\left( s\right) >0$ for all $s>0$ which satisfies%
\begin{equation}
\int_{1}^{\infty }\frac{1}{f\left( t\right) }dt=\infty ,  \label{DY}
\end{equation}%
then the single equation%
\begin{equation}
\left\{ 
\begin{array}{l}
\Delta u=a\left( \left\vert x\right\vert \right) f\left( u\right) \text{ for 
}x\in \mathbb{R}^{N}\text{ (}N\geq 3\text{),} \\ 
\underset{\left\vert x\right\vert \rightarrow \infty }{\lim }u\left(
\left\vert x\right\vert \right) =\infty%
\end{array}%
\right.  \label{dye}
\end{equation}%
has a nonnegative radial solution if and only if $a$ satisfies 
\begin{equation*}
\underset{t\rightarrow \infty }{\lim }\mathcal{A}_{a}\left( t\right) =\infty 
\text{, }\mathcal{A}_{a}\left( t\right)
:=\int_{0}^{t}s^{1-N}\int_{0}^{s}z^{N-1}a(z)dzds.
\end{equation*}%
After a simple computation, we can see that%
\begin{equation*}
\underset{t\rightarrow \infty }{\lim }\mathcal{A}_{a}\left( t\right) =\frac{1%
}{N-2}\int_{0}^{\infty }ra\left( r\right) dr.
\end{equation*}%
However, there is no equivalent results for systems (\ref{11}), where $f_{1}$%
, $f_{2}$ satisfy a condition of the form (\ref{DY}). \ One of the purpose
of this paper is to fill this gap.

Subsequently, Lair \cite{L2} extended the result of \cite{L1} to a more
general case by merely requiring $\alpha \beta \leq 1$, and showed that if $%
\alpha \beta >1$, then (\ref{lair}) has an entire large solution if either (%
\ref{c1l}) and (\ref{c2l}) fails to hold, i.e., $a_{1}$ and $a_{2}$ satisfy
(at least) one of the conditions%
\begin{eqnarray}
\int_{0}^{\infty }ta_{1}\left( t\right) \left(
t^{2-N}\int_{0}^{t}s^{N-3}Q\left( s\right) ds\right) ^{\alpha }dt &<&\infty ,
\label{c3l} \\
\int_{0}^{\infty }ta_{2}\left( t\right) \left(
t^{2-N}\int_{0}^{t}s^{N-3}P\left( s\right) ds\right) ^{\beta }dt &<&\infty .
\label{c4l}
\end{eqnarray}%
To summarises, if $\alpha \beta >1$, a sufficient condition to ensure the
existence of a positive entire large solution for the system (\ref{lair}) is
that $a_{1}$ and $a_{2}$ satisfy (\ref{c3l}) or (\ref{c4l}). Therefore, it
remains unknown whether this is a necessary condition. However, we know from
the reference \cite{DF} that this is not true for the single equation (\ref%
{dye}). The second purpose of this paper is to prove that this does not
happen the systems either.

Finally, we note that if $a_{1}$ and $a_{2}$ satisfy%
\begin{equation}
\begin{array}{llll}
1) & \int_{0}^{\infty }ra_{1}\left( r\right) dr=\infty \text{,} & 2) & 
\int_{0}^{\infty }ra_{2}\left( r\right) dr=\infty ,%
\end{array}%
\text{ }  \label{7l}
\end{equation}%
then they also satisfy both (\ref{c1l}) and (\ref{c2l}), and likewise, if
they satisfy%
\begin{equation}
\begin{array}{cccc}
3) & \int_{0}^{\infty }ra_{1}\left( r\right) dr<\infty \text{,} & 4) & 
\int_{0}^{\infty }ra_{2}\left( r\right) dr<\infty ,%
\end{array}%
\text{ }  \label{8l}
\end{equation}%
then they also satisfy (\ref{c3l}) and (\ref{c4l}). In both cases, however,
the converse is not true. For further results, see for instance, \cite%
{CD,LM,M,NU,NUS} and the references therein.

In the present paper, we are interested in providing a proof to our goals
for a more general class of quasilinear systems of the form (\ref{11}).
This, actually, is the third motivation of our paper since the $\phi _{i}-$%
Laplacian operator appears in mathematical models in nonlinear elasticity,
plasticity, generalized Newtonian fluids, and in quantum physics.

Several results concerning our goals were obtained by Gregorio \cite{DG},
Franchi-Lanconelli and Serrin \cite{FLS}, Hamydy-Massar-Tsouli \cite{HMT},
Keller \cite{K}, Kon'kov \cite{KO}, Jaro\^{s}-Taka\^{s}i \cite{JK},
Losev-Mazepa \cite{LM}, Lieberman \cite{GL}, Li-Zhang-Zhang \cite{Z}, Luthey 
\cite{L}, Mazepa \cite{M}, Naito-Usami \cite{NU,NUS}, Osserman \cite{O},
Smooke \cite{S}, Zhang-Zhou \cite{ZZ} and Zhang \cite{XZ}.

We expect that our work, while currently focussed on a very specific
problem, will lead to general insights and new methods with potential
applications to a much wider class of problems.

Throughout the paper we let $\alpha ,\beta \in \left( 0,\infty \right) $ be
arbitrary parameters. We work under the following assumptions$:$

(A)\quad $a_{1},a_{2}:\left[ 0,\infty \right) \rightarrow \left[ 0,\infty
\right) $ are spherically symmetric continuous functions (i.e.,\textit{\ }$%
a_{i}\left( x\right) =a_{i}\left( \left\vert x\right\vert \right) $ for $%
i=1,2$)\textit{;}

(C1)\quad $f_{1},f_{2}:[0,\infty )\rightarrow \lbrack 0,\infty )$ are
continuous, increasing, $f_{1}\left( 0\right) \cdot f_{2}\left( 0\right)
\geq 0$ and $f_{1}\left( s_{1}\right) \cdot f_{2}\left( s_{2}\right) >0$ for
all $s_{1},s_{2}>0$;

(C2)\quad there exist positive constants $\overline{c}_{1},\overline{c}_{2}$%
, the continuous and increasing functions $g_{1}$, $g_{2}:\left[ 0,\infty
\right) \rightarrow \left[ 0,\infty \right) $ and the continuous functions $%
\overline{\xi }_{1}$, $\overline{\xi }_{2}:\left[ 0,\infty \right)
\rightarrow \left[ 0,\infty \right) $ such that%
\begin{eqnarray}
f_{1}\left( t_{1}\cdot w_{1}\right) &\leq &\overline{c}_{1}g_{1}\left(
t_{1}\right) \cdot \overline{\xi }_{1}\left( w_{1}\right) \text{ }\forall 
\text{ }w_{1}\geq 1\text{ and }\forall \text{ }t_{1}\geq M_{1}\cdot 
\overline{\theta }_{2}(f_{2}\left( \alpha \right) ),  \label{c21} \\
f_{2}\left( t_{2}\cdot w_{2}\right) &\leq &\overline{c}_{2}g_{2}\left(
t_{2}\right) \cdot \overline{\xi }_{2}\left( w_{2}\right) \text{ }\forall 
\text{ }w_{2}\geq 1\text{ and }\forall \text{ }t_{2}\geq M_{2}\cdot 
\overline{\theta }_{1}(f_{1}\left( \beta \right) ),  \label{c22}
\end{eqnarray}%
where $M_{1}\geq \max \left\{ 1,\frac{\beta }{\overline{\theta }%
_{2}(f_{2}\left( \alpha \right) )}\right\} $ and $M_{2}\geq \max \left\{ 1,%
\frac{\alpha }{\overline{\theta }_{1}(f_{1}\left( \beta \right) )}\right\} $;

(C3)\quad there are some constants $\underline{c}_{1},\underline{c}_{2}\in
\left( 0,\infty \right) $ and the continuous functions $\underline{\xi }_{1}$%
, $\underline{\xi }_{2}:\left[ 0,\infty \right) \rightarrow \left[ 0,\infty
\right) $ such that%
\begin{eqnarray}
f_{1}\left( m_{1}w_{1}\right) &\geq &\underline{c}_{1}\underline{\xi }%
_{1}\left( w_{1}\right) \text{ }\forall \text{ }w_{1}\geq 1,  \label{c31} \\
f_{2}\left( m_{2}w_{2}\right) &\geq &\underline{c}_{2}\underline{\xi }%
_{2}\left( w_{2}\right) \text{ }\forall \text{ }w_{2}\geq 1,  \label{c32}
\end{eqnarray}%
where $m_{1}\in \left( 0,\min \left\{ \beta ,\underline{\theta }%
_{2}(f_{2}(\alpha ))\right\} \right) $ and $m_{2}\in \left( 0,\min \left\{
\alpha ,\underline{\theta }_{1}(f_{1}(\beta ))\right\} \right) $.

\section{Main Results}

As announced we start with the formulation of our results. It is convenient
to give some notations needed in the sequel. The reader may just as well
glance through this paper and return to it when necessary 
\begin{eqnarray*}
\overline{\mathcal{A}}_{a_{i}}\left( t\right) &=&\int_{0}^{t}\overline{k}_{i}%
\overline{\psi }_{i}(s^{1-N}\int_{0}^{s}z^{N-1}a_{i}(z)dz)ds\text{, }i=1,2 \\
\overline{P}_{1,2}\left( r\right) &=&\int_{0}^{r}\overline{\psi }_{2}\left( 
\overline{c}_{1}y^{1-N}\int_{0}^{y}t^{N-1}a_{1}(t)\overline{\xi }_{1}\left(
1+\overline{\mathcal{A}}_{a_{2}}\left( t\right) \right) dt\right) dy\text{,}
\\
\overline{P}_{2,1}\left( r\right) &=&\int_{0}^{r}\overline{\psi }_{1}\left( 
\overline{c}_{2}y^{1-N}\int_{0}^{y}t^{N-1}a_{2}(t)\overline{\xi }_{2}\left(
1+\overline{\mathcal{A}}_{a_{1}}\left( t\right) \right) dt\right) dy\text{, }
\\
\text{ }\overline{P}_{1,2}\left( \infty \right) &=&\lim_{r\rightarrow \infty
}\overline{P}_{1,2}\left( r\right) \text{, }\overline{P}_{2,1}\left( \infty
\right) =\lim_{r\rightarrow \infty }\overline{P}_{2,1}\left( r\right) \\
\underline{\mathcal{A}}_{a_{i}}\left( t\right) &=&\int_{0}^{t}\underline{k}%
_{i}\underline{\psi }_{i}(s^{1-N}\int_{0}^{s}z^{N-1}a_{i}(z)dz)ds\text{, }%
i=1,2 \\
\underline{P}_{1,2}\left( r\right) &=&\int_{0}^{r}h_{1}^{-1}\left( 
\underline{c}_{1}y^{1-N}\int_{0}^{y}t^{N-1}a_{1}(t)\underline{\xi }%
_{1}\left( 1+\underline{\mathcal{A}}_{a_{2}}\left( t\right) \right)
dt\right) dy,\text{ } \\
\underline{P}_{2,1}\left( r\right) &=&\int_{0}^{r}h_{2}^{-1}\left( 
\underline{c}_{2}y^{1-N}\int_{0}^{y}t^{N-1}a_{2}(t)\underline{\xi }%
_{2}\left( 1+\underline{\mathcal{A}}_{a_{1}}\left( t\right) \right)
dt\right) dy, \\
\underline{P}_{1,2}\left( \infty \right) &=&\lim_{r\rightarrow \infty }%
\underline{P}_{1,2}\left( r\right) \text{, }\underline{P}_{2,1}\left( \infty
\right) =\lim_{r\rightarrow \infty }\underline{P}_{2,1}\left( r\right) \\
H_{1,2}\left( r\right) &=&\int_{a}^{r}\frac{1}{\overline{\theta }%
_{1}(g_{1}\left( M_{1}\overline{\theta }_{2}(f_{2}\left( t\right) \right) )}%
dt\text{, }H_{1,2}\left( \infty \right) =\lim_{s\rightarrow \infty
}H_{1,2}\left( s\right) \\
H_{2,1}\left( r\right) &=&\int_{b}^{r}\frac{1}{\overline{\theta }%
_{2}(g_{2}\left( M_{2}\overline{\theta }_{1}(f_{1}\left( t\right) \right) )}%
dt\text{, }H_{2,1}\left( \infty \right) =\lim_{s\rightarrow \infty
}H_{2,1}\left( s\right) .
\end{eqnarray*}%
Let us point that 
\begin{equation*}
H_{1,2}^{\prime }(r)=\frac{1}{\overline{\theta }_{1}(g_{1}\left( M_{1}%
\overline{\theta }_{2}(f_{2}\left( r\right) \right) )}>0\text{ for }r>a,
\end{equation*}%
and%
\begin{equation*}
\text{ }H_{2,1}^{\prime }(r)=\frac{1}{\overline{\theta }_{2}(g_{2}\left(
M_{2}\overline{\theta }_{1}(f_{1}\left( r\right) \right) )}>0\text{ for }r>b,
\end{equation*}%
and then $H_{1,2}$ has the inverse function $H_{1,2}^{-1}$ on $%
[0,H_{1,2}(\infty ))$ respectively $H_{2,1}$ has the inverse function $%
H_{2,1}^{-1}$ on $[0,H_{2,1}(\infty ))$.

Having all these notations clearly for the readers, we state the following
first result:

\begin{theorem}
\label{th1}A\textit{ssume that }$H_{1,2}\left( \infty \right) =H_{2,1}\left(
\infty \right) =\infty $ and \textrm{(A),} hold\textit{. Furthermore, if }$%
f_{1}$ and $f_{2}$ satisfy the hypotheses \textrm{(C1)} and \textrm{(C2)}
then the system (\ref{11}) has one positive radial solution 
\begin{equation*}
\left( u,v\right) \in C^{1}\left( \left[ 0,\infty \right) \right) \times
C^{1}\left( \left[ 0,\infty \right) \right) \text{ with }\left( u\left(
0\right) ,v\left( 0\right) \right) =\left( \alpha ,\beta \right) .
\end{equation*}
If in addition, $f_{1}$ and $f_{2}$ satisfy the hypothesis \textrm{(C3)}, $%
\underline{P}_{1,2}\left( \infty \right) =\infty $ and $\underline{P}%
_{2,1}\left( \infty \right) =\infty $ then $\lim_{r\rightarrow \infty
}u\left( r\right) =\infty $ and $\lim_{r\rightarrow \infty }v\left( r\right)
=\infty $. Conversely, if $\underline{\xi }_{i}=\overline{\xi }_{i}$ ($i=1,2$%
), $h_{i}^{-1}=\overline{\psi }_{i}$ ($i=1,2$) and \textrm{(C1), (C2)}, 
\textrm{(C3) }hold true, and $\left( u,v\right) $ is a nonnegative entire
large solution of (\ref{11}) such that $\left( u\left( 0\right) ,v\left(
0\right) \right) =\left( \alpha ,\beta \right) $, then $a_{1}$ and $a_{2}$
satisfy $\underline{P}_{1,2}\left( \infty \right) =\overline{P}_{1,2}\left(
\infty \right) =\infty $ and $\underline{P}_{2,1}\left( \infty \right) =%
\overline{P}_{2,1}\left( \infty \right) =\infty $.
\end{theorem}

Our Theorem \ref{th1} includes all known results about the large solutions
for (\ref{11}) as well as all of the `mixed' cases and therefore gives an
answer for our first goal. Next, we are interested in the existence of
entire bounded radial solutions for the system (\ref{11}).

\begin{theorem}
\label{th12}Suppose\textit{\ that }$H_{1,2}\left( \infty \right)
=H_{2,1}\left( \infty \right) =\infty $ and \textrm{(A),} hold\textit{.
Furthermore, if }$f_{1}$ and $f_{2}$ satisfy the hypotheses \textrm{(C1)}
and \textrm{(C2)} then the system (\ref{11}) has one positive radial
solution 
\begin{equation*}
\left( u,v\right) \in C^{1}\left( \left[ 0,\infty \right) \right) \times
C^{1}\left( \left[ 0,\infty \right) \right) \text{ with }\left( u\left(
0\right) ,v\left( 0\right) \right) =\left( \alpha ,\beta \right) .
\end{equation*}%
Moreover, if $\overline{P}_{1,2}\left( \infty \right) <\infty $ and $%
\overline{P}_{2,1}\left( \infty \right) <\infty $ then $\lim_{r\rightarrow
\infty }u\left( r\right) <\infty $ and $\lim_{r\rightarrow \infty }v\left(
r\right) <\infty $.
\end{theorem}

The next Theorem present the situation when one of the components is bounded
while the other is large.

\begin{theorem}
\label{th13}A\textit{ssume that }$H_{1,2}\left( \infty \right)
=H_{2,1}\left( \infty \right) =\infty $ and \textrm{(A),} hold\textit{.
Furthermore, if }$f_{1}$ and $f_{2}$ satisfy the hypotheses \textrm{(C1)}
and \textrm{(C2)} then the system (\ref{11}) has one positive radial
solution 
\begin{equation*}
\left( u,v\right) \in C^{1}\left( \left[ 0,\infty \right) \right) \times
C^{1}\left( \left[ 0,\infty \right) \right) \text{ with }\left( u\left(
0\right) ,v\left( 0\right) \right) =\left( \alpha ,\beta \right) .
\end{equation*}
Moreover, the following hold:

1.)\quad If in addition, $f_{2}$ satisfy the condition (\ref{c32}), $%
\overline{P}_{1,2}\left( \infty \right) <\infty $ and $\underline{P}%
_{2,1}\left( \infty \right) =\infty $ then $\lim_{r\rightarrow \infty
}u\left( r\right) <\infty $ and $\lim_{r\rightarrow \infty }v\left( r\right)
=\infty .$

2.)\quad If in addition, $f_{1}$ satisfy the condition (\ref{c31}), $%
\underline{P}_{1,2}\left( \infty \right) =\infty $ and $\overline{P}%
_{2,1}\left( \infty \right) <\infty $ then $\lim_{r\rightarrow \infty
}u\left( r\right) =\infty $ and $\lim_{r\rightarrow \infty }v\left( r\right)
<\infty .$
\end{theorem}

We now propose a more refined question concerning the solutions of system (%
\ref{11}).\ In analogy with Theorems \ref{th1}-\ref{th13}, we can also prove
the following three theorems.

\begin{theorem}
\label{th2}Assume that the hypothesis \textrm{(A) }holds. \ If \textrm{(C1),
(C2)}, $\overline{P}_{1,2}\left( \infty \right) <H_{1,2}\left( \infty
\right) <\infty $ and $\overline{P}_{2,1}\left( \infty \right)
<H_{2,1}\left( \infty \right) <\infty $ are satisfied, then the system (\ref%
{11}) has one positive bounded radial solution%
\begin{equation*}
\left( u,v\right) \in C^{1}\left( \left[ 0,\infty \right) \right) \times
C^{1}\left( \left[ 0,\infty \right) \right) \text{ with }\left( u\left(
0\right) ,v\left( 0\right) \right) =\left( \alpha ,\beta \right) .
\end{equation*}
such that 
\begin{equation*}
\left\{ 
\begin{array}{l}
\alpha +\underline{P}_{1,2}\left( r\right) \leq u\left( r\right) \leq
H_{1,2}^{-1}\left( \overline{k}_{1}\overline{P}_{1,2}\left( r\right) \right)
, \\ 
\beta +\underline{P}_{2,1}\left( r\right) \leq v\left( r\right) \leq
H_{2,1}^{-1}\left( \overline{k}_{2}\overline{P}_{2,1}\left( r\right) \right)
.%
\end{array}%
\right. \text{ }
\end{equation*}
\end{theorem}

\begin{theorem}
\label{th21}Assume that the hypothesis \textrm{(A) }holds. \ If \textrm{%
(C1), (C2)}, (\ref{c31}), $H_{1,2}\left( \infty \right) =\infty $, $%
\underline{P}_{1,2}\left( \infty \right) =\infty $ and $\underline{P}%
_{2,1}\left( \infty \right) <H_{2,1}\left( \infty \right) <\infty $ are
satisfied, then the system (\ref{11}) has one positive radial solution 
\begin{equation*}
\left( u,v\right) \in C^{1}\left( \left[ 0,\infty \right) \right) \times
C^{1}\left( \left[ 0,\infty \right) \right) \text{ with }\left( u\left(
0\right) ,v\left( 0\right) \right) =\left( \alpha ,\beta \right) ,
\end{equation*}
such that $\lim_{r\rightarrow \infty }u\left( r\right) =\infty $ and $%
\lim_{r\rightarrow \infty }v\left( r\right) <\infty $.
\end{theorem}

\begin{theorem}
\label{th22}Assume that the hypothesis \textrm{(A) }holds. \ If \textrm{%
(C1), (C2)}, (\ref{c32}), $\underline{P}_{2,1}\left( \infty \right) =\infty $%
, $H_{2,1}\left( \infty \right) =\infty $ and $\overline{P}_{1,2}\left(
\infty \right) <H_{1,2}\left( \infty \right) <\infty $ are satisfied, then
the system (\ref{11}) has one positive radial solution 
\begin{equation*}
\left( u,v\right) \in C^{1}\left( \left[ 0,\infty \right) \right) \times
C^{1}\left( \left[ 0,\infty \right) \right) \text{ with }\left( u\left(
0\right) ,v\left( 0\right) \right) =\left( \alpha ,\beta \right) ,
\end{equation*}%
such that $\lim_{r\rightarrow \infty }u\left( r\right) <\infty $ and $%
\lim_{r\rightarrow \infty }v\left( r\right) =\infty $.
\end{theorem}

\begin{remark}
Our assumptions (O3), (C2) and (C3) are further discussed in the famous book
of Krasnosel'skii and Rutickii \cite{KR} (see also Soria \cite{SO}).
Moreover, the class of nonlinearities considered by Lair \cite{L1}, \cite{L2}
are also included.
\end{remark}

\begin{remark}
(see \cite[Lemma 2.1]{FK}) Suppose $\phi _{i}$ ($i=1,2$) satisfy (O1), (O2)
and

(O4)\quad there exist $l_{i},m_{i}>1$ such that 
\begin{equation*}
l_{i}\leq \frac{\Phi _{i}^{\prime }\left( t\right) \cdot t}{\Phi _{i}\left(
t\right) }\leq m_{i}\text{ for any }t>0\text{, where }\Phi _{i}\left(
t\right) =\int_{0}^{t}\phi _{i}\left( s\right) sds,t>0;
\end{equation*}

(O5)\quad there exist $a_{0}^{i}$, $a_{1}^{i}>0$ such that 
\begin{equation*}
a_{0}^{i}\leq \frac{\Phi _{i}^{\prime \prime }\left( t\right) \cdot t}{\Phi
_{i}^{\prime }\left( t\right) }\leq a_{1}^{i}\text{ for any }t>0\text{.}
\end{equation*}%
Then, the assumption (\ref{ineq}) holds with%
\begin{equation*}
\overline{\psi }_{i}=\underline{\psi }_{i}=h_{i}^{-1}\text{, }\underline{k}%
_{i}=\overline{k}_{i}=1\text{, }\underline{\theta }_{i}(t)=\min \left\{
t^{1/m_{i}},t^{1/l_{i}}\right\} \text{, }\overline{\theta }_{i}(t)=\max
\left\{ t^{1/m_{i}},t^{1/l_{i}}\right\} .
\end{equation*}
\end{remark}

We would like to point that:

\begin{remark}[see \protect\cite{FK} for more information]
The function $\Phi _{i}$ appears in a lot of physical applications, such as:

Nonlinear Elasticity: $\Phi _{i}\left( t\right) =\left( 1+t^{2}\right)
^{p}-1 $, $\phi _{i}\left( t\right) =2p\left( 1+t^{2}\right) ^{p-1}$, $t>0$
and $p>\frac{1}{2}$;

Plasticity: $\Phi _{i}\left( t\right) =t^{p}\left( \ln \left( 1+t\right)
\right) ^{q}$, $\phi _{i}\left( t\right) =\frac{\ln ^{q-1}\left( t+1\right) 
}{t+1}\left[ \left( pt^{p-1}+qt^{q-2}\right) \ln \left( t+1\right) +qt^{p-1}%
\right] $, $t>0$, $p>1$ and $q>0$;

Generalized Newtonian fluids: $\Phi _{i}\left( t\right)
=\int_{0}^{t}s^{1-p}\left( \sinh ^{-1}s\right) ^{q}ds$, $\phi _{i}\left(
t\right) =\allowbreak t^{-p}\func{arcsinh}^{q}t$, $t>0$, $0\leq p\leq 1$ and 
$q>0$;

Plasma Physics: $\Phi _{i}\left( t\right) =\frac{t^{p}}{p}+\frac{t^{q}}{q}$, 
$\phi _{i}\left( t\right) =t^{p-2}+t^{q-2}$ where $t>0$ and $1<p<q$.
\end{remark}

\begin{remark}
Let%
\begin{equation*}
M_{1}^{+}=\underset{t\in \left[ 0,\infty \right) }{\sup }\int_{0}^{t}%
\overline{k}_{i}\overline{\psi }_{2}(s^{1-N}\int_{0}^{s}z^{N-1}a_{2}(z)dz)ds
\end{equation*}%
and%
\begin{equation*}
M_{2}^{+}=\underset{t\in \left[ 0,\infty \right) }{\sup }\int_{0}^{t}%
\overline{k}_{1}\overline{\psi }_{1}(s^{1-N}\int_{0}^{s}z^{N-1}a_{1}(z)dz)ds.
\end{equation*}%
The following situations improve our theorems:

a)\quad If $M_{1}^{+}\in \left( 0,\infty \right) $ then the condition (\ref%
{c21}) is not necessary but $H_{1,2}\left( r\right) $ must be replaced by 
\begin{equation}
H_{1,2}\left( r\right) =\int_{a}^{r}\frac{1}{\overline{\theta }%
_{1}(f_{1}\left( M_{1}\overline{\theta }_{2}(f_{2}\left( t\right) \right) )}%
dt\text{, }M_{1}\geq \max \left\{ 1,\frac{\beta }{\overline{\theta }%
_{2}(f_{2}\left( \alpha \right) )}\right\} \cdot \left( 1+M_{1}^{+}\right) ,
\label{t1}
\end{equation}%
and therefore $\overline{P}_{1,2}\left( r\right) =\int_{0}^{r}\overline{\psi 
}_{2}\left( \overline{c}_{1}y^{1-N}\int_{0}^{y}t^{N-1}a_{1}(t)dt\right) dy$.

b)\quad If $M_{2}^{+}\in \left( 0,\infty \right) $ then the condition (\ref%
{c22}) is not necessary but $H_{2,1}\left( r\right) $ must be replaced by 
\begin{equation}
H_{2,1}\left( r\right) =\int_{b}^{r}\frac{1}{\overline{\theta }%
_{2}(f_{2}\left( M_{2}\overline{\theta }_{1}(f_{1}\left( t\right) \right) )}%
dt\text{, }M_{2}\geq \max \left\{ 1,\frac{\alpha }{\overline{\theta }%
_{1}(f_{1}\left( \beta \right) )}\right\} \cdot \left( 1+M_{2}^{+}\right) .
\label{t12}
\end{equation}%
and therefore $\overline{P}_{2,1}\left( r\right) =\int_{0}^{r}\overline{\psi 
}_{1}\left( \overline{c}_{2}y^{1-N}\int_{0}^{y}t^{N-1}a_{2}(t)dt\right) dy$.

c)\quad If $M_{1}^{+}\in \left( 0,\infty \right) $ and $M_{2}^{+}\in \left(
0,\infty \right) $ then the conditions (\ref{c21}) \ and (\ref{c22}) are not
necessary but $H_{1,2}\left( r\right) $ and $H_{2,1}\left( r\right) $ must
be replaced by (\ref{t1}) and (\ref{t12}). Here $\overline{P}_{1,2}\left(
r\right) $ and $\overline{P}_{2,1}\left( r\right) $ are defined as in a), b).

d)\quad If $m_{1}\geq 1$ then $\underline{c}_{1}=1$ and $\underline{\xi }%
_{1}=f_{1}$.

e)\quad If $m_{2}\geq 1$ then $\underline{c}_{2}=1$ and $\underline{\xi }%
_{2}=f_{2}$.

f)\quad If $m_{1}\geq 1$ and $m_{2}\geq 1$ then $\underline{c}_{1}=%
\underline{c}_{2}=1$, $\underline{\xi }_{1}=f_{1}$ and $\underline{\xi }%
_{2}=f_{2}$.
\end{remark}

\section{Proof of the main results}

The first important tool in our proof is a variant of the Arzel\`{a}--Ascoli
Theorem.

\subsection{The Arzel\`{a}--Ascoli Theorem}

Let $r_{1},r_{2}\in \mathbb{R}$ with $r_{1}\leq r_{2}$ and 
\begin{equation*}
\left( K=\left[ r_{1},r_{2}\right] ,d_{K}\left( x,y\right) \right)
\end{equation*}
be a compact metric space, with the metric $d_{K}\left( x,y\right)
=\left\vert x-y\right\vert $, and let 
\begin{equation*}
C\left( \left[ r_{1},r_{2}\right] \right) =\left\{ g:\left[ r_{1},r_{2}%
\right] \rightarrow \mathbb{R}\left\vert g\text{ is continuous on }\left[
r_{1},r_{2}\right] \right. \right\}
\end{equation*}%
denote the space of real valued continuous functions on $\left[ r_{1},r_{2}%
\right] $ and for any $g\in C\left( \left[ r_{1},r_{2}\right] \right) $, let 
\begin{equation*}
\left\Vert g\right\Vert _{\infty }=\underset{x\in \left[ r_{1},r_{2}\right] }%
{\max }\left\vert g\left( x\right) \right\vert
\end{equation*}%
be the maximum norm on $C\left( \left[ r_{1},r_{2}\right] \right) $.

\begin{remark}
Let $g^{1},g^{2}\in C\left( \left[ r_{1},r_{2}\right] \right) $. If $d\left(
g^{1},g^{2}\right) =\left\Vert g^{1}-g^{2}\right\Vert _{\infty }$ then $%
\left( C\left( \left[ r_{1},r_{2}\right] \right) ,d\right) $ is a complete
metric space.
\end{remark}

\begin{definition}
We say that the sequence $\left\{ g_{n}\right\} _{n\in \mathbb{N}}$ from $%
C\left( \left[ r_{1},r_{2}\right] \right) $ is bounded if there exists a
positive constant $C<\infty $ such that $\left\Vert g_{n}\left( x\right)
\right\Vert _{\infty }\leq C$ for each $x\in \left[ r_{1},r_{2}\right] $.
(Equivalently: $\left\vert g_{n}\left( x\right) \right\vert \leq C$ for each 
$x\in \left[ r_{1},r_{2}\right] $ and $n\in \mathbb{N}^{\ast }$).
\end{definition}

\begin{definition}
We say that the sequence $\left\{ g_{n}\right\} _{n\in \mathbb{N}}$ from $%
C\left( \left[ r_{1},r_{2}\right] \right) $ is equicontinuous if for any
given $\varepsilon >0$, there exists a number $\delta >0$ (which depends
only on $\varepsilon $) such that 
\begin{equation*}
\left\vert g_{n}\left( x\right) -g_{n}\left( y\right) \right\vert
<\varepsilon \text{ for all }n\in \mathbb{N}
\end{equation*}%
whenever $d_{K}\left( x,y\right) <\delta $ for every $x,y\in \left[
r_{1},r_{2}\right] $.
\end{definition}

\begin{definition}
Let $\left\{ g_{n}\right\} _{n\in \mathbb{N}}$ be a family of functions
defined on $\left[ r_{1},r_{2}\right] $. The sequence $\left\{ g_{n}\right\}
_{n\in \mathbb{N}}$ converges uniformly to $g\left( x\right) $ if for every $%
\varepsilon >0$ there is an $N$ (which depends only on $\varepsilon $) such
that%
\begin{equation*}
\left\vert g_{n}\left( x\right) -g\left( x\right) \right\vert <\varepsilon 
\text{ for all }n>N\text{ and }x\in \left[ r_{1},r_{2}\right] .
\end{equation*}
\end{definition}

\begin{theorem}[Arzel\`{a}--Ascoli theorem]
\label{arzela}If a sequence $\left\{ g_{n}\right\} _{n\in \mathbb{N}}$ in $%
C\left( \left[ r_{1},r_{2}\right] \right) $ is bounded and equicontinuous
then it has a subsequence $\left\{ g_{n_{k}}\right\} _{k\in \mathbb{N}}$
which converges uniformly to $g\left( x\right) $ on $C\left( \left[
r_{1},r_{2}\right] \right) $.
\end{theorem}

\subsection{Proof of Theorems \protect\ref{th1}- \protect\ref{th13}}

Radially symmetric solutions of the problem (\ref{11}) correspond to
solutions of the ordinary differential equations system%
\begin{equation}
\left\{ 
\begin{array}{c}
\left( r^{N-1}\phi _{1}\left( \left\vert u^{\prime }\left( r\right)
\right\vert \right) u^{\prime }\left( r\right) \right) ^{\prime
}=r^{N-1}a_{1}(r)f_{1}(v\left( r\right) )\text{ on }\left[ 0,\infty \right) 
\text{,} \\ 
\left( r^{N-1}\phi _{2}\left( \left\vert v^{\prime }\left( r\right)
\right\vert \right) v^{\prime }\left( r\right) \right) ^{\prime
}=r^{N-1}a_{2}(r)f_{2}(u\left( r\right) )\text{ on }\left[ 0,\infty \right) 
\text{,}%
\end{array}%
\right.  \label{ords}
\end{equation}%
subject to the initial conditions $\left( u\left( 0\right) ,v\left( 0\right)
\right) =\left( \alpha ,\beta \right) $ and $\left( u^{\prime }\left(
0\right) ,v^{\prime }\left( 0\right) \right) =\left( 0,0\right) $, since $%
\left( u\left( r\right) ,v\left( r\right) \right) $ is a radially symmetric
positive entire solution of the system (\ref{11}). Integrating (\ref{ords})
from $0$ to $r$, we obtain%
\begin{equation}
\left\{ 
\begin{array}{c}
\phi _{1}(\left\vert u^{\prime }(r)\right\vert )u^{\prime }(r)=\frac{1}{%
r^{N-1}}\int_{0}^{r}r^{N-1}a_{1}\left( r\right) f_{1}\left( v\left( s\right)
\right) ds\text{, on }\left[ 0,\infty \right) , \\ 
\phi _{2}(\left\vert v^{\prime }(r)\right\vert )v^{\prime }(r)=\frac{1}{%
r^{N-1}}\int_{0}^{r}r^{N-1}a_{2}\left( r\right) f_{2}\left( u\left( s\right)
\right) ds\text{, on }\left[ 0,\infty \right) .%
\end{array}%
\right.  \label{ssord}
\end{equation}%
Taking into account the equations (\ref{ssord}), it is easy to see that $%
u\left( r\right) $ is an increasing function on $\left[ 0,\infty \right) $
of the radial variable $r$, and the same conclusion holds for $v\left(
r\right) $. Thus, for radial solutions of the system (\ref{ords}) we seek
for solutions of the system of integral equations 
\begin{equation}
\left\{ 
\begin{array}{l}
u(r)=\alpha
+\int_{0}^{r}h_{1}^{-1}(t^{1-N}\int_{0}^{t}s^{N-1}a_{1}(s)f_{1}(v(s))ds)dt,%
\quad r\geq 0, \\ 
v(r)=\beta
+\int_{0}^{r}h_{2}^{-1}(t^{1-N}\int_{0}^{t}s^{N-1}a_{2}(s)f_{2}(u(s))ds)dt,%
\quad r\geq 0.%
\end{array}%
\right.  \label{ords2}
\end{equation}%
The system (\ref{ords2}) can be solved by using successive approximation. We
define inductively $\{u_{m}\}_{m\geq 0}$ and $\{v_{m}\}_{m\geq 0}$ on $\left[
0,\infty \right) $ as follows 
\begin{equation}
\left\{ 
\begin{array}{l}
u_{0}(r)=\alpha \text{, }v_{0}(r)=\beta , \\ 
u_{m}(r)=\alpha
+\int_{0}^{r}h_{1}^{-1}(t^{1-N}%
\int_{0}^{t}s^{N-1}a_{1}(s)f_{1}(v_{m-1}(s))ds)dt,\quad r\geq 0, \\ 
v_{m}(r)=\beta
+\int_{0}^{r}h_{2}^{-1}(t^{1-N}%
\int_{0}^{t}s^{N-1}a_{2}(s)f_{2}(u_{m}(s))ds)dt,\quad r\geq 0.%
\end{array}%
\right.  \label{rec}
\end{equation}%
Obviously, for all $r\geq 0$ and $m\in {\mathbb{N}}$ it holds that $%
u_{m}(r)\geq \alpha $, $v_{m}(r)\geq \beta $ and $v_{0}\leq v_{1}$. Our
assumptions yield $u_{1}(r)\leq u_{2}(r)$, for all $r\geq 0$, so $%
v_{1}(r)\leq v_{2}(r)$, for all $r\geq 0$. Continuing on this line of
reasoning, we obtain that the sequences $\{u_{m}\}_{m}$ and $\{v_{m}\}_{m}$
are increasing on $[0,\infty )$.

We next establish bounds for the non-decreasing sequences $\{u_{m}\}_{m}$
and $\{v_{m}\}_{m}$. \ From (\ref{rec}) we obtain the following inequalities 
\begin{eqnarray}
v_{m}(r) &=&\beta
+\int_{0}^{r}h_{2}^{-1}(t^{1-N}%
\int_{0}^{t}s^{N-1}a_{2}(s)f_{1}(u_{m}(s))ds)dt  \notag \\
&\leq &\beta
+\int_{0}^{r}h_{2}^{-1}(f_{2}(u_{m}(t))t^{1-N}%
\int_{0}^{t}z^{N-1}a_{2}(z)dz)dt  \notag \\
&\leq &\beta +\int_{0}^{r}\overline{k}_{2}\overline{\theta }_{2}(f_{2}\left(
u_{m}(t)\right) \overline{\psi }_{2}(t^{1-N}\int_{0}^{t}z^{N-1}a_{2}(z)dz)dt
\notag \\
&\leq &\beta +\overline{\theta }_{2}(f_{2}\left( u_{m}(r)\right)
)\int_{0}^{r}\overline{k}_{2}\overline{\psi }_{2}(t^{1-N}%
\int_{0}^{t}z^{N-1}a_{2}(z)dz)dt  \label{inex1} \\
&\leq &\overline{\theta }_{2}(f_{2}\left( u_{m}(r)\right) (\frac{\beta }{%
\overline{\theta }_{2}(f_{2}\left( u_{m}(r)\right) }+\overline{\mathcal{A}}%
_{a_{2}}\left( r\right) )  \notag \\
&\leq &\overline{\theta }_{2}(f_{2}\left( u_{m}(r)\right) (\frac{\beta }{%
\overline{\theta }_{2}(f_{2}\left( \alpha \right) )}+\overline{\mathcal{A}}%
_{a_{2}}\left( r\right) )  \notag \\
&\leq &M_{1}\overline{\theta }_{2}(f_{2}\left( u_{m}(r)\right) (1+\overline{%
\mathcal{A}}_{a_{2}}\left( r\right) )  \notag
\end{eqnarray}%
and, in the same vein%
\begin{eqnarray}
u_{m}(r) &=&\alpha
+\int_{0}^{r}h_{1}^{-1}(t^{1-N}%
\int_{0}^{t}s^{N-1}a_{1}(s)f_{1}(v_{m-1}(s))ds)dt  \notag \\
&\leq &\alpha
+\int_{0}^{r}h_{1}^{-1}(t^{1-N}%
\int_{0}^{t}s^{N-1}a_{1}(s)f_{1}(v_{m}(s))ds)dt  \label{inex2} \\
&\leq &M_{2}\overline{\theta }_{1}(f_{1}\left( v_{m}(r)\right) (1+\overline{%
\mathcal{A}}_{a_{1}}\left( r\right) ).  \notag
\end{eqnarray}%
Moreover, using (\ref{inex1}), by an elementary computation it follows that 
\begin{align}
u_{m}^{\prime }(r)& \leq h_{1}^{-1}\left(
r^{1-N}\int_{0}^{r}s^{N-1}a_{1}(s)f_{1}(v_{m}(s))ds\right)  \notag \\
& \leq h_{1}^{-1}\left( r^{1-N}\int_{0}^{r}s^{N-1}a_{1}(s)f_{1}\left( M_{1}%
\overline{\theta }_{2}(f_{2}\left( u_{m}(s)\right) (1+\overline{\mathcal{A}}%
_{a_{2}}\left( s\right) )\right) ds\right)  \notag \\
& \leq h_{1}^{-1}\left( r^{1-N}\int_{0}^{r}s^{N-1}a_{1}(s)\overline{c}%
_{1}g_{1}\left( M_{1}\overline{\theta }_{2}(f_{2}\left( u_{m}(s)\right)
\right) \overline{\xi }_{1}\left( 1+\overline{\mathcal{A}}_{a_{2}}\left(
s\right) \right) ds\right)  \label{exin} \\
& \leq h_{1}^{-1}\left( g_{1}\left( M_{1}\overline{\theta }_{2}(f_{2}\left(
u_{m}(r)\right) \right) \overline{c}_{1}r^{1-N}\int_{0}^{r}s^{N-1}a_{1}(s)%
\overline{\xi }_{1}\left( 1+\overline{\mathcal{A}}_{a_{2}}\left( s\right)
\right) ds\right)  \notag \\
& \leq \overline{k}_{1}\overline{\theta }_{1}(g_{1}\left( M_{1}\overline{%
\theta }_{2}(f_{2}\left( u_{m}(r)\right) \right) )\overline{\psi }_{2}\left( 
\overline{c}_{1}r^{1-N}\int_{0}^{r}s^{N-1}a_{1}(s)\overline{\xi }_{1}\left(
1+\overline{\mathcal{A}}_{a_{2}}\left( s\right) \right) ds\right) .  \notag
\end{align}%
Arguing as above, but now with the second inequality (\ref{inex2}), one can
show that 
\begin{eqnarray}
v_{m}^{\prime }(r) &=&h_{2}^{-1}\left(
r^{1-N}\int_{0}^{r}s^{N-1}a_{2}(s)f_{1}(u_{m-1}(s))ds\right)  \label{exin2}
\\
&\leq &\overline{k}_{2}\overline{\theta }_{2}(g_{2}\left( M_{2}\overline{%
\theta }_{1}(f_{1}\left( v_{m}(r)\right) \right) )\overline{\psi }_{1}\left( 
\overline{c}_{2}r^{1-N}\int_{0}^{r}s^{N-1}a_{2}(s)\overline{\xi }_{2}\left(
1+\overline{\mathcal{A}}_{a_{1}}\left( s\right) \right) ds\right) .  \notag
\end{eqnarray}%
Combining the previous relations (\ref{exin}) and (\ref{exin2}), we further
obtain 
\begin{eqnarray}
\frac{\left( u_{1}^{m}\left( r\right) \right) ^{\prime }}{\overline{\theta }%
_{1}(g_{1}\left( M_{1}\overline{\theta }_{2}(f_{2}\left( u_{m}(r)\right)
\right) )} &\leq &\overline{k}_{1}\overline{\psi }_{2}\left( \overline{c}%
_{1}r^{1-N}\int_{0}^{r}s^{N-1}a_{1}(s)\overline{\xi }_{1}\left( 1+\overline{%
\mathcal{A}}_{a_{2}}\left( s\right) \right) ds\right) ,  \label{mat1} \\
\frac{\left( u_{2}^{m}\left( r\right) \right) ^{\prime }}{\overline{\theta }%
_{2}(g_{2}\left( M_{2}\overline{\theta }_{1}(f_{1}\left( v_{m}(r)\right)
\right) )} &\leq &\overline{k}_{2}\overline{\psi }_{1}\left( \overline{c}%
_{2}r^{1-N}\int_{0}^{r}s^{N-1}a_{2}(s)\overline{\xi }_{2}\left( 1+\overline{%
\mathcal{A}}_{a_{1}}\left( s\right) \right) ds\right) .  \label{mat2}
\end{eqnarray}%
Integrating the inequalities (\ref{mat1}) and (\ref{mat2}) from $0$ to $r$,
yields that 
\begin{equation}
\int_{a}^{u_{m}\left( r\right) }\frac{\overline{k}_{1}^{-1}}{\overline{%
\theta }_{1}(g_{1}\left( M_{1}\overline{\theta }_{2}(f_{2}\left( t\right)
\right) )}dt\leq \overline{P}_{1,2}\left( r\right) \text{, }%
\int_{b}^{v_{m}\left( r\right) }\frac{\overline{k}_{2}^{-1}}{\overline{%
\theta }_{2}(g_{2}\left( M_{2}\overline{\theta }_{1}(f_{1}\left( t\right)
\right) )}dt\leq \overline{P}_{2,1}\left( r\right) .  \label{genord}
\end{equation}%
Also, going back to the setting of $H_{1,2}$ and $H_{2,1}$ we rewrite (\ref%
{genord}) as 
\begin{equation}
H_{1,2}\left( u_{m}(r)\right) \leq \overline{k}_{1}\overline{P}_{1,2}\left(
r\right) \text{ and }H_{2,1}\left( v_{m}(r)\right) \leq \overline{k}_{2}%
\overline{P}_{2,1}\left( r\right) ,  \label{ints}
\end{equation}%
which plays a basic role in the proof of our main results. Since $H_{1,2}$
(resp. $H_{2,1}$) is a bijection with the inverse function $H_{1,2}^{-1}$
(resp. $H_{2,1}^{-1}$) strictly increasing on $\left[ 0,\infty \right) $,
the inequalities (\ref{ints}) can be reformulated as 
\begin{equation}
u_{m}(r)\leq H_{1,2}^{-1}\left( \overline{k}_{1}\overline{P}_{1,2}\left(
r\right) \right) \text{ and }v_{m}(r)\leq H_{2,1}^{-1}\left( \overline{k}_{2}%
\overline{P}_{2,1}\left( r\right) \right) .  \label{int}
\end{equation}%
So, we have found upper bounds for $\left\{ u_{m}(r)\right\} _{m\geq 1}$ and 
$\left\{ v_{m}(r)\right\} _{m\geq 1}$which are dependent of $r$. We point to
the reader that the corresponding estimates (\ref{int}) are sometimes
essential.

Next we prove that the sequences $\left\{ u_{m}(r)\right\} _{m\geq 1}$ and $%
\left\{ v_{m}(r)\right\} _{m\geq 1}$ are bounded and equicontinuous on $%
\left[ 0,c_{0}\right] $ for arbitrary $c_{0}>0$. To do this, we take 
\begin{equation*}
C_{1}=H_{1,2}^{-1}\left( \overline{k}_{1}\overline{P}_{1,2}\left(
c_{0}\right) \right) \text{ and }C_{2}=H_{2,1}^{-1}\left( \overline{k}_{2}%
\overline{P}_{2,1}\left( c_{0}\right) \right)
\end{equation*}%
and since $\left( u_{m}(r)\right) ^{^{\prime }}\geq 0$ and $\left(
v_{m}\left( r\right) \right) ^{^{\prime }}\geq 0$ it follows that 
\begin{equation*}
u_{m}(r)\leq u_{m}\left( c_{0}\right) \leq C_{1}\text{ and }v_{m}\left(
r\right) \leq v_{m}\left( c_{0}\right) \leq C_{2}.
\end{equation*}%
We have proved that $\left\{ u_{m}(r)\right\} _{m\geq 1}$ and $\left\{
v_{m}(r)\right\} _{m\geq 1}$ are bounded on $\left[ 0,c_{0}\right] $ for
arbitrary $c_{0}>0$. Using this fact in (\ref{exin}) and (\ref{exin2}) we
show that the same is true for $\left( u_{m}(r)\right) ^{\prime }$ and $%
\left( v_{m}(r)\right) ^{\prime }$. By construction we verify that 
\begin{eqnarray}
u_{m}^{\prime }(r) &=&h_{1}^{-1}\left(
r^{1-N}\int_{0}^{r}s^{N-1}a_{1}(s)f_{1}(v_{m-1}(s))ds\right)  \notag \\
&\leq &h_{1}^{-1}\left(
r^{1-N}\int_{0}^{r}s^{N-1}a_{1}(s)f_{1}(v_{m}(s))ds\right)  \notag \\
&\leq &h_{1}^{-1}\left( \int_{0}^{r}a_{1}(s)f_{1}(v_{m-1}(s))ds\right) 
\notag \\
&\leq &h_{1}^{-1}\left( \left\Vert a_{1}\right\Vert _{\infty
}\int_{0}^{r}f_{1}(v_{m-1}(s))ds\right)  \label{maj} \\
&\leq &h_{1}^{-1}\left( \left\Vert a_{1}\right\Vert _{\infty
}f_{1}(C_{2})\int_{0}^{r}ds\right)  \notag \\
&\leq &h_{1}^{-1}\left( \left\Vert a_{1}\right\Vert _{\infty
}f_{1}(C_{2})c_{0}\right) \text{ on }\left[ 0,c_{0}\right] .  \notag
\end{eqnarray}%
We follow the argument used in (\ref{maj}) to obtain 
\begin{equation*}
\left( v_{m}(r)\right) ^{\prime }\leq h_{2}^{-1}\left( \left\Vert
a_{2}\right\Vert _{\infty }f_{2}(C_{1})c_{0}\right) \text{ on }\left[ 0,c_{0}%
\right] .
\end{equation*}%
Summarizing, we have found that 
\begin{eqnarray*}
\left( u_{1}^{m}\left( r\right) \right) ^{\prime } &\leq &h_{1}^{-1}\left(
\left\Vert a_{1}\right\Vert _{\infty }f_{1}(C_{2})c_{0}\right) \text{ on }%
\left[ 0,c_{0}\right] , \\
\left( v_{m}(r)\right) ^{\prime } &\leq &h_{2}^{-1}\left( \left\Vert
a_{2}\right\Vert _{\infty }f_{2}(C_{1})c_{0}\right) \text{ on }\left[ 0,c_{0}%
\right] .
\end{eqnarray*}%
Finally, it remains to prove that $\left\{ u_{m}(r)\right\} _{m\geq 1}$ and $%
\left\{ v_{m}(r)\right\} _{m\geq 1}$ are equicontinuous on $\left[ 0,c_{0}%
\right] $ for arbitrary $c_{0}>0$. Let $\varepsilon _{1}$, $\varepsilon
_{2}>0$ be arbitrary. To verify equicontinuity on $\left[ 0,c_{0}\right] $
observe that the mean value theorem yields 
\begin{eqnarray*}
\left\vert u_{m}\left( x\right) -u_{m}\left( y\right) \right\vert
&=&\left\vert \left( u_{m}\left( \zeta _{1}\right) \right) ^{\prime
}\right\vert \left\vert x-y\right\vert \leq h_{1}^{-1}\left( \left\Vert
a_{1}\right\Vert _{\infty }f_{1}(C_{2})c_{0}\right) \left\vert
x-y\right\vert , \\
\left\vert v_{m}\left( x\right) -v_{m}\left( y\right) \right\vert
&=&\left\vert \left( v_{m}\left( \zeta _{2}\right) \right) ^{\prime
}\right\vert \left\vert x-y\right\vert \leq h_{2}^{-1}\left( \left\Vert
a_{2}\right\Vert _{\infty }f_{2}(C_{1})c_{0}\right) \left\vert
x-y\right\vert ,
\end{eqnarray*}%
for all $n\in \mathbb{N}$ and all $x,y\in \left[ 0,c_{0}\right] $ and for
some $\zeta _{1}$, $\zeta _{2}$. Then it suffices to take 
\begin{equation*}
\delta _{1}=\frac{\varepsilon _{1}}{h_{1}^{-1}\left( \left\Vert
a_{1}\right\Vert _{\infty }f_{1}(C_{2})c_{0}\right) }\text{ and }\delta _{2}=%
\frac{\varepsilon _{2}}{h_{2}^{-1}\left( \left\Vert a_{2}\right\Vert
_{\infty }f_{2}(C_{1})c_{0}\right) }
\end{equation*}%
to see that $\left\{ u_{m}(r)\right\} _{m\geq 1}$ and $\left\{
v_{m}(r)\right\} _{m\geq 1}$ are equicontinuous on $\left[ 0,c_{0}\right] $.

Since $\left\{ u_{m}(r)\right\} _{m\geq 1}$ and $\left\{ v_{m}(r)\right\}
_{m\geq 1}$ are bounded and equicontinuous on $\left[ 0,c_{0}\right] $ we
can apply the Arzel\`{a}--Ascoli theorem with $\left[ r_{1},r_{2}\right] =%
\left[ 0,c_{0}\right] $. Thus, there exists a subsequence, denoted $\left\{
\left( u_{m^{1}}(r),v_{\overline{m}^{1}}(r)\right) \right\} $ that converges
uniformly on $\left[ 0,1\right] \times \left[ 0,1\right] $. Let%
\begin{equation*}
\left( u_{m^{1}}\left( r\right) ,v_{\overline{m}^{1}}\left( r\right) \right) 
\overset{\left( m^{1},\overline{m}^{1}\right) \rightarrow \infty }{%
\rightarrow }\left( u_{1}\left( r\right) ,v_{1}\left( r\right) \right) \text{
uniformly on }\left[ 0,1\right] .\text{ }
\end{equation*}%
Likewise, the subsequence $\left\{ \left( u_{m^{1}}\left( r\right) ,v_{%
\overline{m}^{1}}\left( r\right) \right) \right\} $ is bounded and
equicontinuous on the interval $\left[ 0,2\right] $. Hence, it must contain
a convergent subsequence 
\begin{equation*}
\left\{ \left( u_{m^{2}}\left( r\right) ,v_{\overline{m}^{2}}\left( r\right)
\right) \right\} ,
\end{equation*}
that converges uniformly on $\left[ 0,2\right] \times \left[ 0,2\right] $.
Let%
\begin{equation*}
\left( u_{m^{2}}\left( r\right) ,v_{\overline{m}^{2}}\left( r\right) \right) 
\overset{\left( m^{2},\overline{m}^{2}\right) \rightarrow \infty }{%
\rightarrow }\left( u_{2}\left( r\right) ,v_{2}\left( r\right) \right) \text{
uniformly on }\left[ 0,2\right] \times \left[ 0,2\right] .
\end{equation*}%
Note that 
\begin{equation*}
\left\{ u_{m^{2}}\left( r\right) \right\} \subseteq \left\{ u_{m^{1}}\left(
r\right) \right\} \subseteq \left\{ u_{m}\left( r\right) \right\} ^{m\geq 2}%
\text{ and }\left\{ v_{\overline{m}^{2}}\left( r\right) \right\} \subseteq
\left\{ v_{\overline{m}^{1}}\left( r\right) \right\} \subseteq \left\{
v_{m}\left( r\right) \right\} ^{m\geq 2}.
\end{equation*}%
These imply 
\begin{equation*}
u_{2}\left( r\right) =u_{1}\left( r\right) \text{ and }v_{2}\left( r\right)
=v_{1}\left( r\right) \text{ on }\left[ 0,1\right] .
\end{equation*}%
Proceeding in this fashion we obtain a countable collection of subsequences
such that 
\begin{equation*}
\left\{ u_{m^{n}}\right\} \subseteq ....\subseteq \left\{ u_{m^{1}}\left(
r\right) \right\} \subseteq \left\{ u_{m}\left( r\right) \right\} _{m\geq n}%
\text{ }
\end{equation*}%
and%
\begin{equation*}
\left\{ v_{\overline{m}^{n}}\right\} \subseteq ....\subseteq \left\{ v_{%
\overline{m}^{1}}\left( r\right) \right\} \subseteq \left\{ v_{m}\left(
r\right) \right\} _{m\geq n}
\end{equation*}%
and a sequence $\left\{ \left( u_{n}\left( r\right) ,v_{n}\left( r\right)
\right) \right\} $ such that%
\begin{equation*}
\text{ }%
\begin{array}{lll}
\left( u_{n}\left( r\right) ,v_{n}\left( r\right) \right) \in C\left[ 0,n%
\right] \times C\left[ 0,n\right] & \text{for } & n=1,2,3,... \\ 
\left( u_{n}\left( r\right) ,v_{n}\left( r\right) \right) =\left(
u_{1}\left( r\right) ,v_{1}\left( r\right) \right) & \text{for } & r\in 
\left[ 0,1\right] \\ 
\left( u_{n}\left( r\right) ,v_{n}\left( r\right) \right) =\left(
u_{2}\left( r\right) ,v_{2}\left( r\right) \right) & \text{for } & r\in 
\left[ 0,2\right] \\ 
... & \text{...} & ... \\ 
\left( u_{n}\left( r\right) ,v_{n}\left( r\right) \right) =\left(
u_{n-1}\left( r\right) ,v_{n-1}\left( r\right) \right) & \text{for } & r\in 
\left[ 0,n-1\right] .%
\end{array}%
\end{equation*}%
Together, these observations show that there exists a sequence $\left\{
\left( u_{n}\left( r\right) ,v_{n}\left( r\right) \right) \right\} $ that
converges to $\left( u\left( r\right) ,v\left( r\right) \right) $ on $\left[
0,\infty \right) $ satisfying 
\begin{equation*}
\left( u_{n}\left( r\right) ,v_{n}\left( r\right) \right) =\left( u\left(
r\right) ,v\left( r\right) \right) \text{ if }0\leq r\leq n.
\end{equation*}%
This convergence is uniformly on bounded intervals, implying $\left( u\left(
r\right) ,v\left( r\right) \right) \in C\left[ 0,\infty \right) \times C%
\left[ 0,\infty \right) $, and moreover, the family $\left\{ \left(
u_{n}\left( r\right) ,v_{n}\left( r\right) \right) \right\} $ is also
equicontinuous. The solution $\left( u\left( r\right) ,v\left( r\right)
\right) $ constructed in this way is radially symmetric.

Going back to the system (\ref{ords}), the radial solutions of (\ref{11})
are the solutions of the ordinary differential equations system (\ref{ords}%
). We conclude that radial solutions of (\ref{11}) with $u\left( 0\right)
=\alpha ,$ $v\left( 0\right) =\beta $ satisfy:%
\begin{eqnarray}
u(r) &=&\alpha
+\int_{0}^{r}h_{1}^{-1}(t^{1-N}\int_{0}^{t}s^{N-1}a_{1}(s)f_{1}(v(s))ds)dt,r%
\geq 0,  \label{eq1} \\
v(r) &=&\beta
+\int_{0}^{r}h_{2}^{-1}(t^{1-N}\int_{0}^{t}s^{N-1}a_{2}(s)f_{2}(u(s))ds)dt,%
\text{ }r\geq 0.  \label{eq2}
\end{eqnarray}%
We are now ready to give a complete proof of the Theorems \ref{th1}-\ref%
{th13}.

\subsubsection{\textbf{Proof of Theorem \protect\ref{th1} completed:} \ }

\ From (\ref{eq2}) we obtain the following inequalities%
\begin{eqnarray*}
v(r) &=&\beta
+\int_{0}^{r}h_{2}^{-1}(t^{1-N}\int_{0}^{t}s^{N-1}a_{2}(s)f_{2}(u(s))ds)dt \\
&\geq &\beta +\int_{0}^{r}h_{2}^{-1}(f_{2}(\alpha
)z^{1-N}\int_{0}^{z}s^{N-1}a_{2}(s)ds)dz \\
&\geq &\beta +\underline{\theta }_{2}(f_{2}(\alpha ))\underline{\mathcal{A}}%
_{a_{2}}\left( r\right) \\
&\geq &m_{1}(1+\underline{\mathcal{A}}_{a_{2}}\left( r\right) ),
\end{eqnarray*}%
and, in the same vein%
\begin{eqnarray*}
u(r) &=&\alpha
+\int_{0}^{r}h_{1}^{-1}(t^{1-N}\int_{0}^{t}s^{N-1}a_{1}(s)f_{1}(v(s))ds)dt \\
&\geq &m_{2}(1+\underline{\mathcal{A}}_{a_{1}}\left( r\right) ).
\end{eqnarray*}%
If $\underline{P}_{1,2}\left( \infty \right) =\underline{P}_{2,1}\left(
\infty \right) =\infty $, we observe that 
\begin{eqnarray}
u\left( r\right) &=&\alpha
+\int_{0}^{r}h_{1}^{-1}(t^{1-N}\int_{0}^{t}s^{N-1}a_{1}(s)f_{1}(v(s))ds)dt 
\notag \\
&\geq &\alpha +\int_{0}^{r}h_{1}^{-1}\left(
y^{1-N}\int_{0}^{y}t^{N-1}a_{1}(t)f_{1}\left( m_{1}(1+\underline{\mathcal{A}}%
_{a_{2}}\left( t\right) \right) dt\right) dy  \notag \\
&\geq &\alpha +\int_{0}^{r}h_{1}^{-1}\left( \underline{c}_{1}y^{1-N}%
\int_{0}^{y}t^{N-1}a_{1}(t)\underline{\xi }_{1}\left( 1+\underline{\mathcal{A%
}}_{a_{2}}\left( t\right) \right) dt\right) dy  \label{i11} \\
&\geq &\alpha +\int_{0}^{r}h_{1}^{-1}\left( \underline{c}_{1}y^{1-N}%
\int_{0}^{y}t^{N-1}a_{1}(t)\underline{\xi }_{1}\left( 1+\underline{\mathcal{A%
}}_{a_{2}}\left( t\right) \right) dt\right) dy  \notag \\
&=&\alpha +\underline{P}_{1,2}\left( r\right) .  \notag
\end{eqnarray}%
Analogously, we refine the strategy above to prove: \ 
\begin{eqnarray*}
v\left( r\right) &\geq &\beta +\int_{0}^{r}h_{2}^{-1}\left( \underline{c}%
_{2}y^{1-N}\int_{0}^{y}t^{N-1}a_{2}(t)\underline{\xi }_{2}\left( 1+%
\underline{\mathcal{A}}_{a_{1}}\left( t\right) \right) dt\right) dy \\
&=&\beta +\underline{P}_{2,1}\left( r\right) ,
\end{eqnarray*}%
and passing to the limit as $r\rightarrow \infty $ in (\ref{i11}) and in the
above inequality we conclude that%
\begin{equation*}
\lim_{r\rightarrow \infty }u\left( r\right) =\lim_{r\rightarrow \infty
}v\left( r\right) =\infty ,
\end{equation*}%
which yields the result. In order to prove the converse let $\left(
u,v\right) $ be an entire large radial solution of (\ref{11}) such that $%
\left( u,v\right) =\left( \alpha ,\beta \right) $. Then, $\left( u,v\right) $
satisfy%
\begin{eqnarray*}
u(r) &=&\alpha
+\int_{0}^{r}h_{1}^{-1}(t^{1-N}\int_{0}^{t}s^{N-1}a_{1}(s)f_{1}(v(s))ds)dt,%
\text{ }r\geq 0, \\
v(r) &=&\beta
+\int_{0}^{r}h_{2}^{-1}(t^{1-N}\int_{0}^{t}s^{N-1}a_{2}(s)f_{2}(u(s))ds)dt,%
\text{ }r\geq 0,
\end{eqnarray*}%
and, so%
\begin{equation}
\text{ }H_{1,2}\left( u\left( r\right) \right) \leq \overline{k}_{1}%
\overline{P}_{1,2}\left( r\right) \text{ and }H_{2,1}\left( v\left( r\right)
\right) \leq \overline{k}_{2}\overline{P}_{2,1}\left( r\right) .
\label{conv}
\end{equation}%
By passing to the limit as $r\rightarrow \infty $ in (\ref{conv}) we find
that $a_{1}$ and $a_{2}$ satisfy $\overline{P}_{1,2}\left( \infty \right) =%
\overline{P}_{2,1}\left( \infty \right) =\infty $, since $\left( u,v\right) $
is large and $H_{1,2}\left( \infty \right) =H_{2,1}\left( \infty \right)
=\infty $. This completes the proof. We next consider:

\subsubsection{\textbf{Proof of Theorem \protect\ref{th12} completed:} \ }

If $\overline{P}_{1,2}\left( \infty \right) <\infty $ and $\overline{P}%
_{2,1}\left( \infty \right) <\infty ,$ then using the same arguments as in (%
\ref{eq1}) and (\ref{eq2}) we can see that%
\begin{equation*}
u\left( r\right) \leq H_{1,2}^{-1}\left( \overline{k}_{1}\overline{P}%
_{1,2}\left( \infty \right) \right) <\infty \text{ and }v\left( r\right)
\leq H_{2,1}^{-1}\left( \overline{k}_{2}\overline{P}_{2,1}\left( \infty
\right) \right) <\infty \text{ for all }r\geq 0.
\end{equation*}%
Hence $\left( u,v\right) $ is bounded\textbf{\ }and\textbf{\ }this completes
the proof.

\subsubsection{\textbf{Proof of Theorem \protect\ref{th13} completed:}}

\textbf{Case 1):} By an analysis similar to the Theorems \ref{th1} and \ref%
{th13} above, we have that 
\begin{equation*}
u\left( r\right) \leq H_{1,2}^{-1}\left( \overline{k}_{1}\overline{P}%
_{1,2}\left( \infty \right) \right) <\infty \text{ and }v\left( r\right)
\geq b+\overline{k}_{2}\underline{P}_{2,1}\left( r\right) .
\end{equation*}%
So, if 
\begin{equation*}
\overline{P}_{1,2}\left( \infty \right) <\infty \text{ and }\underline{P}%
_{2,1}\left( \infty \right) =\infty
\end{equation*}%
then 
\begin{equation*}
\lim_{r\rightarrow \infty }u\left( r\right) <\infty \text{ and }%
\lim_{r\rightarrow \infty }v\left( r\right) =\infty .
\end{equation*}%
In order, to complete the proofs it remains to proceed to the

\textbf{Case 2): }In this case, we invoke the proof of Theorem \ref{th12}.
An easy computation yields that 
\begin{equation}
u\left( r\right) \geq \alpha +\overline{k}_{1}\underline{P}_{1,2}\left(
r\right) \text{ and }v\left( r\right) \leq H_{2,1}^{-1}\left( \overline{k}%
_{2}\overline{P}_{2,1}\left( r\right) \right) .  \label{t2}
\end{equation}%
Our conclusion follows by letting $r\rightarrow \infty $ in (\ref{t2}).

\subsection{Proof of Theorems \protect\ref{th2}- \protect\ref{th22}}

\subsubsection{\textbf{Proof of Theorem \protect\ref{th2} completed: }}

We deduce from (\ref{ints}) and the conditions of the theorem that 
\begin{equation*}
\begin{array}{l}
H_{1,2}\left( u_{m}\left( r\right) \right) \leq \overline{k}_{1}\overline{P}%
_{1,2}\left( \infty \right) <\overline{k}_{1}H_{1,2}\left( \infty \right)
<\infty , \\ 
H_{2,1}\left( v_{m}\left( r\right) \right) \leq \overline{k}_{2}\overline{P}%
_{2,1}\left( \infty \right) <\overline{k}_{2}H_{2,1}\left( \infty \right)
<\infty .%
\end{array}%
\end{equation*}%
On the other hand, since $H_{1,2}^{-1}$ and $H_{2,1}^{-1}$ are strictly
increasing on $\left[ 0,\infty \right) $, we find out that%
\begin{equation*}
u_{m}\left( r\right) \leq H_{1,2}^{-1}\left( \overline{k}_{1}\overline{P}%
_{1,2}\left( \infty \right) \right) <\infty \text{ and }v_{m}\left( r\right)
\leq H_{2,1}^{-1}\left( \overline{k}_{2}\overline{P}_{2,1}\left( \infty
\right) \right) <\infty ,
\end{equation*}%
and then the non-decreasing sequences $\left\{ u_{m}\left( r\right) \right\}
_{m\geq 1}$ and $\left\{ v_{m}\left( r\right) \right\} _{m\geq 1}$ are
bounded above for all $r\geq 0$ and all $m$. Putting these two facts
together yields 
\begin{equation*}
\left( u_{m}\left( r\right) ,v_{m}\left( r\right) \right) \rightarrow \left(
u\left( r\right) ,v\left( r\right) \right) \text{ as }m\rightarrow \infty
\end{equation*}%
and the limit functions $u$ and $v$ are positive entire bounded radial
solutions of system (\ref{11}).\textbf{\ }

\subsubsection{\textbf{Proof of Theorem \protect\ref{th21} and \protect\ref%
{th22} completed: }}

It is a straightforward adaptation of the above proofs.

\bigskip

\end{document}